%% file: BM+SRW_covering3.tex
\documentclass[amsppt,12pt]{amsart}
\usepackage{curves}
\usepackage{amsmath}
\newtheorem {theorem}{Theorem}[section]
\newtheorem {definition}{Definition}[section]

\newtheorem {lemma}[theorem]{Lemma}
\newtheorem {corollary}[theorem]{Corollary}

\newcounter{conjecture}
\newcounter{remark}

\newcommand{\eqnsection}{
    \renewcommand{\theequation}{\thesection.\arabic{equation}}
    \makeatletter
    \csname @addtoreset\endcsname{equation}{section}
    \makeatother}

\include{macros}

\newcommand{\ls}[1]
    {\dimen0=\fontdimen6\the\font \lineskip=#1\dimen0
\advance\lineskip.5\fontdimen5\the\font \advance\lineskip-\dimen0
\lineskiplimit=.9\lineskip \baselineskip=\lineskip
\advance\baselineskip\dimen0 \normallineskip\lineskip
\normallineskiplimit\lineskiplimit \normalbaselineskip\baselineskip
\ignorespaces }

\begin{document}

\bibliographystyle{amsplain}

\title[An LIL for cover times of disks]
{An LIL for cover times of disks \\
by planar random walk and Wiener sausage}

\author[J. Ben Hough\, \,and \,\, Yuval Peres]
{J. Ben Hough$^*$\,\, and\,\, Yuval Peres$^*$}

\date{October 18, 2004.
\newline\indent
$^*$The authors gratefully
acknowledge the financial support from NSF grants $\#$DMS-0104073 and $\#$DMS-0244479}

\begin{abstract}
\noindent
Let $R_n$ be the radius of the largest disk covered after $n$ steps of a simple random walk.  We prove that almost surely $\limsup_{n \rightarrow \infty}(\log R_n)^2/(\log n \log_3 n) = 1/4$, where $\log_3$ denotes 3 iterations of the $\log$ function.  This is motivated by a question of Erd\H{o}s and Taylor.  We also obtain the analogous result for the Wiener sausage, refining a result of Meyre and Werner.   
\end{abstract}

\maketitle

\section{Introduction}

In this paper, we consider a planar simple random walk starting from the origin and determine a sharp asymptotic upper bound for the growth of the radius of the largest discrete disk centered at the origin that is covered after $n$ steps of the walk.  More precisely, let $D_r = D(0,r) \cap \Z^2$ be the discrete disk of radius $r$, $S(n)$ a planar simple random walk starting from the origin, and $R_n = \sup \{r:D_r \subset S[0,n] \}$ the radius of the largest discrete disk centered at the origin covered by time $n$.  It is well known that the planar simple random walk is recurrent so $R_n \rightarrow \infty$ as $n \rightarrow \infty$.  We prove that almost surely
\begin{equation} \label{eqn1}
	\limsup_{n \rightarrow \infty} \frac{(\log R_n)^2}{\log n \log_3 n} = \frac{1}{4},
\end{equation}
where $\log_k$ denotes $k$ iterations of the $\log$ function.  We also consider the analogous covering problem for the planar Wiener sausage of radius 1, and prove that the same asymptotic result holds in this case as well.  Specifically, if $\BB_t$ denotes a planar Brownian motion starting from the origin, the Wiener sausage of radius 1 up to time $t$ is defined to be 
\begin{equation}
	\SS_t = \bigcup_{s\in[0,t]} D(\BB_s,1). \nonumber
\end{equation}
We define $\RR_t = \sup\{r:D(0,r) \subset \SS_t \}$ and prove that almost surely
\begin{equation}
	\limsup_{t \rightarrow \infty} \frac{(\log \RR_t)^2}{\log t \log_3 t} = \frac{1}{4}.
\end{equation}

A problem related to the asymptotic behavior of $R_n$ was posed as early as 1960 by Erd\H{o}s and Taylor.  In \cite[p.153]{erdos} they ask: \textquotedblleft How quickly does the function $f(n)$ need to increase so that in an infinite plane random walk, with probability 1, all the lattice points within a distance $n$ of the origin will be entered by the walk before $f(n)$ steps except for finitely many values of $n$?"  The asymptotic growth of $\RR_t$ has been studied by Meyre and Werner in \cite{werner}, who establish the existence of $c \in [1/8,1]$ so that almost surely $\limsup_{t \rightarrow \infty} \frac{(\log \RR_t)^2}{\log t \log_3 t} = c$.  

The present paper relies heavily on the use of excursions between large concentric circles to control covering times.  This technique has been used effectively in many recent works on the subject, see \cite{DPRZ}, \cite{lawler} and \cite{werner} for example.  Although it is seemingly difficult to calculate the cover time for a disk directly, one may say very precisely how many excursions will occur between large concentric circles prior to the covering time of the disk.  The cover time may then be determined indirectly by studying the expected duration of the requisite number of excursions.  We derive the results above by establishing sufficiently tight bounds on the tails of the relevant distributions.  Much of the work lies in estimating the number of excursions necessary to cover a disk.  Lemma \ref{RWest} gives this estimate for the simple random walk, and is derived using results from \cite{DPRZ} and \cite{lawler}.  Lemma \ref{WSest} gives an analogous estimate for the Wiener sausage, and is derived using results from \cite{DPRZ}.

This article is divided into two main sections.  In the first section we study the planar simple random walk process and prove (\ref{eqn1}).  The second section demonstrates that analogous results hold for the planar Wiener sausage.

\section{Asymptotic results for the random walk}
We begin with some definitions.  Recall that $S(n)$ is a simple random walk in the plane started at zero, and $D_r = D(0,r) \cap \Z^2$ is the disk of radius $r$ in $\Z^2$.  Define random variables
\begin{eqnarray}
R_n = \sup \{ r:D_r \subset S[0,n] \} \nonumber \\
T_r = \inf \{n:D_r \subset S[0,n] \} \nonumber
\end{eqnarray}
so that $R_n$ is the largest disk centered at zero that is covered by time $n$, and $T_r$ is the time necessary to cover $D_r$.  The following discussion will demonstrate that
\begin{equation}
\limsup_{n \rightarrow \infty} \frac{(\log R_n)^2}{\log n \log_3 n} = \frac{1}{4} \; \mathrm{a.s.}
\end{equation}
The proof has two main components, we first check that $\frac{1}{4}$ is an upper bound for the $\limsup$, and then deduce that it is a lower bound as well.  The proofs of a few technical estimates will be saved for the end.

\subsection{Establishing that $\frac{1}{4}$ is an upper bound}

For notational simplicity, define
\begin{equation}
  f(x) = \exp{\{[\lambda\log x \log_3 x]^{1/2}\}},
\end{equation}
$\wp(x) = x (\log x)^3$ and $t_n = e^{\alpha^n}$ where $\alpha>1$ is arbitrary.  It is sufficient to show that whenever $\lambda>\frac{1}{4}$ there exists $\alpha = \alpha(\lambda) >1$ so that $\P\{A_n \; \mathrm{i.o.}\}=0$ where $A_n$ is the event that $R_{t_{n+1}} \geq f(t_n)$.  We compute probabilities of events determined by $A_n$'s by counting excursions.  Let $\partial D_r$ be the boundary of $D_r$ in $\Z^2$:
\begin{equation}
  \partial D_r = \{z \in D_r^c \cap \Z^2 : |y-z| = 1 \; \mathrm{for \;some} \; y \in D_r\}.
\end{equation}
\begin{definition}
We say that an {\em excursion} from $\partial D_r$ to $\partial D_p$ occurs between times $s<t$ if $S(s) \in \partial D_r$, $S(t) \in \partial D_p$ and the times $s$ and $t$ satisfy:
\begin{itemize}
\item[1.]  $t = \min\{n:n>s \; {\rm and} \; S(n) \in \partial D_p \}$
\item[2.]  $s = \min\{n:n>\tau \; {\rm and} \; S(n) \in \partial D_r \}$ where $\tau = \max\{n:n<t \; {\rm and } \;S(n) \in \partial D_p \}$ (we take the max of the empty set to be 0).
\end{itemize}
\end{definition}
The key result linking excursions to our covering problem is that if $N_r$ is the number of excursions between $\partial D_{2r}$ and $\partial D_{\wp(r)}$ that occur before $T_r$, and $\phi_r = \frac{(\log r)^2}{\log_2 r}$ then
\begin{equation} \label{eqn6}
\frac{N_r}{\phi_r} \rightarrow \frac{2}{3}
\end{equation}
in probability as $r \rightarrow \infty$ (\cite[Lemma 5.1]{DPRZ} and \cite[Lemma 2.4, Theorem 2.5]{lawler}).

With this result in mind, we are inspired to introduce the events $\tilde{A}_n$ that at least $(\frac{2}{3} - \epsilon_1)\phi_{f(t_n)}$ excursions from $\partial D_{2 f(t_n)}$ to $\partial D_{\wp(f(t_n))}$ occur before the walk first hits $\partial D_{t_{n+1}^{1/2 + \epsilon_2}}$.  The events $A_n$ and $\tilde{A}_n$ should then be comparable since a random walk requires roughly $r^2$ steps to first hit $\partial D_r$.  We now quote \cite[Exercise 1.6.8]{lawler2} which will allow us to compute the probabilities of the events $\tilde{A}_n$.  In what follows, $\zeta(r)$ gives the hitting time of $\partial D_r$.

\begin{lemma} \label{lem22}
If $\rho<r<P$ and $x \in \partial D_r$,
\begin{equation}
\P^x\{ \zeta(\rho) < \zeta(P) \} = \frac{\log(P) - \log(r) + O(1/\rho)}{\log(P) - \log(\rho)},
\end{equation}
where the $O(\cdot)$ term is bounded uniformly in $x$, $\rho$ and $P$.
\end{lemma}
The following corollary follows easily, see also \cite[Proposition 1.6.7]{lawler2}
\begin{corollary} \label{cor23}
If $0<r<P$ and $x \in \partial D_r$ then: \label{simple}
\begin{equation}
\P^x\{\zeta(0)<\zeta(P)\} \geq \frac{1}{16}\left(\frac{\log(P) - \log(r) +O(1)}{\log(P)}\right).
\end{equation}
\end{corollary}

From the above lemma, we see immediately that

\begin{eqnarray}\label{eqn9}
\P(\tilde{A}_n) &=& \left[ \frac{\log t_{n+1}^{1/2+\epsilon_2} - \log(\wp(f(t_n)))] + O(1)}{\log  t_{n+1}^{1/2+\epsilon_2} - \log [2f(t_n)]} \right]^{\lfloor (2/3 - \epsilon_1) \phi_{f(t_n)} \rfloor} \nonumber \\
&=& \left[1 - \frac{\log [\wp(f(t_n))]-\log[2 f(t_n)] + O(1)}{\log t_{n+1}^{1/2 + \epsilon_2} - \log[2 f(t_n)]} \right]^{\lfloor (2/3 - \epsilon_1)\phi_{f(t_n)} \rfloor } \nonumber \\
&=&\left[1 - \frac{3\log_2[f(t_n)] + O(1)}{\log t_{n+1}^{1/2 + \epsilon_2} - \log[2 f(t_n)]} \right]^{\lfloor (2/3 - \epsilon_1)\phi_{f(t_n)} \rfloor }.\label{ugly}
\end{eqnarray}
We use the following two estimates to simplify the expression for $\P(\tilde{A}_n)$.

\begin{lemma} \label{est1}
For $|\epsilon_n|<\frac{1}{2}|a_n|$ and $\delta_n<\frac{1}{2}|b_n|$ we have 
\begin{equation}
\frac{a_n-\epsilon_n}{b_n+\delta_n} = \frac{a_n}{b_n} + O\left(\frac{\epsilon_n}{b_n}\right) + O\left(\frac{a_n\delta_n}{b_n^2}\right). \nonumber
\end{equation}
\end{lemma}
\begin{lemma}\label{est2}
If $\alpha_n^2\beta_n$ and $\beta_n\epsilon_n$ are bounded, $\alpha_n,\frac{\epsilon_n}{\alpha_n} \rightarrow 0$ and $| \delta_n | <1$ then:
\begin{equation}
[1-\alpha_n + \epsilon_n]^{\beta_n + \delta_n} = e^{-\alpha_n \beta_n}[1 + O(\alpha_n^2\beta_n) + O(\beta_n \epsilon_n) + O(\alpha_n \delta_n)].\nonumber
\end{equation}
\end{lemma}

Using Lemma \ref{est1} to simplify the fraction in (\ref{ugly}) and then applying Lemma \ref{est2}, we may simplify the expression for $\P(\tilde{A}_n)$:
\begin{eqnarray}
\P(\tilde{A}_n) &=& \exp \left[ \frac{-4(1-3\epsilon_1/2)}{1+2\epsilon_2}\frac{(\log f(t_n))^2}{\log t_{n+1}} \right] \left\{1 + O\left(\frac{[\log_2 f(t_n)]^2}{(\log t_{n+1})^2}\phi_{f(t_n)}\right) + O\left(\frac{\phi_{f(t_n)}}{\log(t_{n+1})}\right) \right. \nonumber \\
&&+ \left. O \left( \frac{\log_2 f(t_n)}{\log t_{n+1}} \right) \right\}.
\end{eqnarray}

It is easy to verify that
\begin{equation}
\lim_{n \rightarrow \infty} \frac{[\log_2 f(t_n)]^2}{(\log t_{n+1})^2}\phi_{f(t_n)}= 0, \; \mathrm{and} \; \lim_{n \rightarrow \infty} \frac{\phi_{f(t_n)}}{\log(t_{n+1})} = 0, \; \mathrm{and}\; \lim_{n \rightarrow \infty} \frac{\log_2 f(t_n)}{\log t_{n+1}} =0
\end{equation}
so $\sum_{n=1}^\infty \P(\tilde{A}_n)$ converges if and only if

\begin{equation} \label{eqn13}
\sum_{n=1}^\infty{ \exp \left[ \frac{-4(1-3\epsilon_1/2)}{1+2\epsilon_2}\frac{(\log f(t_n))^2}{\log t_{n+1}} \right]}
\end{equation}
converges.  Now compute:
\begin{eqnarray}
\sum_{n=1}^\infty{ \exp \left[ \frac{-4(1-3\epsilon_1/2)}{1+2\epsilon_2}\frac{(\log f(t_n))^2}{\log t_{n+1}} \right]} &=&
\sum_{n=1}^\infty{\exp \left[ \frac{-4(1-3\epsilon_1/2)}{1+2\epsilon_2}\frac{\lambda\alpha^n \log(n\log \alpha)}{\alpha^{n+1}} \right]} \nonumber \\
&=& C \sum_{n=1}^\infty{n^{\frac{-4\lambda}{\alpha}\frac{1-3\epsilon_1/2}{1+2\epsilon_2}}}.
\end{eqnarray}
It follows that if $\lambda>\frac{1}{4}$ we may choose $\epsilon_1>0$ and $\epsilon_2>0$ and $\alpha>1$ so that the above sum converges, and hence $\P\left\{\tilde{A}_n \; \mathrm{i.o.}\right\} = 0$.

To finish the argument, we must demonstrate that the event $\{A_n \; \mathrm{i.o.}\} - \{\tilde{A}_n \; \mathrm{i.o.}\}$ has probability zero.  Observe that $\{A_n \; \mathrm{i.o.}\} - \{\tilde{A}_n \; \mathrm{i.o.}\} \subset C_1 \cup D_1$ where
\begin{eqnarray}
C_1 &=& \left\{\zeta(t_{n}^{1/2 + \epsilon_2})<t_{n}\; \mathrm{i.o.}\right\} \nonumber \\
D_1 &=& \left\{ N_{f(t_n)} < \left(2/3 - \epsilon_1\right)\phi_{f(t_n)} \; \mathrm{i.o.} \right\}. \nonumber
\end{eqnarray}

The fact that $C_1$ has probability zero follows from the estimate \cite[p.146]{lawler2}:
\begin{lemma} \label{lem26}
$\P\left\{\zeta(r) \notin \left(r^{2 - \epsilon}, r^{2 + \epsilon}\right)\right\} = O(e^{-r^\gamma})$, for some $\gamma = \gamma(\epsilon)>0$.
\end{lemma}

And the following lemma shows that $D_1$ is also null.
\begin{lemma}\label{RWest}
$\P\left\{\frac{N_r}{\phi_r} \notin \left(\frac{2}{3} - \epsilon, \frac{2}{3} + \epsilon\right)\right\} \leq O\left(\frac{1}{(\log_2 r)^2}\right)$.
\end{lemma}

Indeed the claim follows from the following computation:
\begin{equation} 
\sum_{n=1}^\infty{\frac{1}{(\log_2 f(t_n))^2}} = \sum_{n=1}^\infty{\frac{4}{(n\log{\alpha} + \log \lambda + \log_3 \alpha^n)^2}}<\infty
\end{equation}
and an application of the Borel-Cantelli lemma.  The proof of this lemma is rather involved, and will be given at the end of the section.

\subsection{Establishing that $\frac{1}{4}$ is a lower bound}

This proof is nearly identical to the one just given.  We define the events $H_n$ and $\tilde{H}_n$ so that
$H_n$ is the event that $R_{t_n}>f(t_n)$ and 
$\tilde{H}_n$ is the event that there are at least $(2/3 + \epsilon_1)\phi_{f(t_n)}$ excursions between $\partial D_{2f(t_n)}$ and $\partial D_{\wp(f(t_n))}$ before $\zeta(t_n^{1/2-\epsilon_2})$.  The task now is to show that if $\lambda<\frac{1}{4}$ then $\P\{H_n \; \mathrm{i.o.}\} = 1$.  We will do this by first showing that we may choose $\epsilon_1>0$ and $\epsilon_2>0$ so that $\P\{\tilde{H}_n \; \mathrm{i.o.}\} =1$ and then checking that for such $\epsilon_1$ and $\epsilon_2$ the two sets are identical up to null sets.

To this end, we define sigma fields $\F_n = \sigma \{ S(k): k<\zeta(t_n^{1/2 - \epsilon_2}) \}$, and check that if $\lambda < \frac{1}{4}$ then for all sufficiently small $\epsilon_1>0$ and $\epsilon_2>0$, we have $\sum_{n=1}^\infty{\P(\tilde{H}_n | \F_{n-1})} = \infty$ with probability one.  Observe that by the strong Markov property:
\begin{eqnarray}\label{eqn16}
  \P(\tilde{H}_n|\F_{n-1}) &\geq& \min_{x \in \partial D_{t_{n-1}^{1/2 - \epsilon_2}}} \P^x(\zeta(0)<\zeta(t_n^{1/2 - \epsilon_2})) \P(\tilde{H}_n) \nonumber \\
 &\geq& \frac{1}{16}\left(\frac{\log[t_n^{1/2 - \epsilon_2}] - \log[t_{n-1}^{1/2 - \epsilon_2}] +O(1)}{\log[t_n^{1/2 - \epsilon_2}]}\right) \P(\tilde{H}_n) \nonumber \\
&\geq&\frac{1}{32}\left(1-\frac{1}{\alpha}\right)\P(\tilde{H}_n),
\end{eqnarray}
where the last inequality follows from Corollary \ref{simple}.

Now compute:
\begin{eqnarray}\label{eqn17}
\sum_{n=1}^\infty \P(\tilde{H}_n) &=& \sum_{n=1}^\infty{\left[\frac{\log(t_n^{1/2-\epsilon_2}) - \log(\wp(f(t_n))) + O(1)}{\log(t_n^{1/2 - \epsilon_2}) - \log(2f(t_n))}\right]^{\lfloor (2/3 + \epsilon_1)\phi_{f(t_n)}\rfloor}} \nonumber \\
&=& \sum_{n=1}^\infty{\left[1-\frac{3\log_2(f(t_n)) + O(1)}{(1/2-\epsilon_2)\log(t_n) - \log(2f(t_n))}\right]^{\lfloor (2/3 + \epsilon_1)\phi_{f(t_n)}\rfloor} } \nonumber \\
&=& \sum_{n=1}^\infty \exp \left[\frac{-4(1 + 3\epsilon_1/2)}{1-2\epsilon_2}\frac{[\log f(t_n)]^2}{\log t_n}\right] \left\{ 1 + O\left(\frac{[\log_2 f(t_n)]^2}{(\log t_n )^2}\phi_{f(t_n)}\right) + O\left(\frac{\phi_{f(t_n)}}{\log t_n}\right) \right. \nonumber\\
& & + \left. O\left(\frac{\log_2 f(t_n)}{\log t_n}\right) \right\}.
\end{eqnarray}
Both error terms tend to zero as $n \to \infty$, so it suffices to consider the sum
\begin{eqnarray} \label{eqn18}
\sum_{n=1}^\infty{\exp \left[\frac{-4(1 + 3\epsilon_1/2)}{1-2\epsilon_2}\frac{[\log f(t_n)]^2}{\log t_n}\right]} &=& \sum_{n=1}^\infty{\exp \left[\frac{-4(1 + 3\epsilon_1/2)}{1-2\epsilon_2}\frac{\lambda\alpha^n\log_2 \alpha^n}{\alpha^n} \right]} \nonumber \\
&=& C \sum_{n=1}^\infty{n^{-4\lambda \frac{1+3\epsilon_1/2}{1-2\epsilon_2}}}.
\end{eqnarray}
If $\lambda<\frac{1}{4}$, we may choose $\epsilon_1>0$ and $\epsilon_2>0$ small enough to that the above sum diverges and thus the sum  $\sum_{n=1}^\infty \P(\tilde{H}_n|\F_{n-1})$ diverges a.s.  Now applying the Borel-Cantelli lemma for filtrations \cite[Corollary 7.20]{K} we may conclude that $\P\{\tilde{H}_n \; \mathrm{i.o.}\} = 1$.

To complete the proof, we must show that $\P( \{\tilde{H}_n \; \mathrm{i.o.} \} - \{H_n \; \mathrm{i.o} \}) = 0$.  Observe that$\{\tilde{H}_n \; \mathrm{i.o.} \} - \{ H_n \; \mathrm{i.o} \} \subset C_2 \cup D_2$, where
\begin{eqnarray}
C_2 &=& \{ \zeta(t_n^{1/2 - \epsilon_2}) > t_n \; \mathrm{i.o.} \} \nonumber \\
D_2 &=& \{ N_{f(t_n)} \geq (2/3 + \epsilon_1)\phi_{f(t_n)} \; \mathrm{i.o.} \}. \nonumber
\end{eqnarray}
Lemmas \ref{lem26} and \ref{RWest} show that both these sets are null, so the proof is complete.

\subsection{Proof of Lemma \ref{RWest}}

The method of proof is to consider the simpler problem of covering $D_r$ by random sets that are covered during $\mathrm{i.i.d.}$ excursions and then to show that our problem does not differ substantially from this situation.  More precisely, define $k(0,r) = 0$, $s(0,r) = \zeta(\wp(r))$ and for $j \geq 1$:
\begin{eqnarray}
	k(i,r) &=& \inf \{j>s(i-1,r):S(j) \in \partial D_{2r} \} \nonumber \\
	s(i,r) &=& \inf \{j>k(i,r):S(j) \in \partial D_{\wp(r)} \}. \nonumber
\end{eqnarray}
Let $A(j,r) = S[k(j,r),s(j,r)]\cap D_r$ consist of the points in $D_r$ covered by the $j^{\mathrm{th}}$ excursion of the random walk.  Take $C(j,r)$, $E(j,r)$ and $F(j,r)$ to be $\mathrm{i.i.d.}$ subsets of $D_r$ chosen according to the following distribution:  pick $x \in \partial D_{2r}$ according to harmonic measure, start a random walk at $x$, stop the walk when it reaches $\partial D_{\wp(r)}$ and choose the points in $D_r$ that are visited by the walk.  Harmonic measure, $H_r$, on $\partial D_r$ is the hitting probability from infinity:
\begin{equation}
	H_r(x) = \lim_{|y| \to \infty} \P^y\{S(\zeta(r)) = x\}. \nonumber
\end{equation}
Let $U(k,r) = \bigcup_{j = 0}^k A(j,r)$.  Also, define $V_C(k,r) = \bigcup_{j=0}^k C(j,r)$ and $V_E$, $V_F$ similarly.  In what follows, we frequenty conserve notation by writing $V_i(x,r)$ when we should more properly write $V_i(\lfloor x \rfloor,r)$.  We now prove an estimate for the $\mathrm{i.i.d.}$ covering process:

\begin{lemma} \label{estimate}
$\P\left\{V_i((\frac{2}{3} - \epsilon)\phi_r,r) = D_r\right\} \leq O\left(\frac{1}{\phi_r}\right)$ and $\P\left\{V_i((\frac{2}{3} + \epsilon)\phi_r,r) \neq D_r\right\} \leq O\left(\frac{1}{\phi_r}\right)$, where $i = C,E$ or $F$.
\end{lemma}
\begin{proof}
From \cite[Proposition 2.6]{lawler} we have the following result:
\begin{lemma}
	Suppose for some $0<u,c,\alpha<\infty$ and all $r$ 
\begin{equation} \label{eqn19}
\P\{V_i(u\phi_r,r) \neq D_r\} \geq c(\log r)^{-\alpha}.  
\end{equation}
Then for all $v<u$, 
\begin{equation}
\lim_{r \rightarrow \infty} \P\{V_i(v\phi_r,r) = D_r \} = 0. 
\end{equation}
\end{lemma}
Applying the proof given in \cite{lawler} to a subsequence of $r$'s, we deduce that if equation (\ref{eqn19}) holds for infinitely many $r$'s then for all $v<u$, 
\begin{equation}
\liminf_{r \rightarrow \infty} \P\{V_i(v \phi_r,r) = D_r \} = 0. \label{extension}
\end{equation}
From (\ref{eqn6}) we know that $\frac{N_r}{\phi_r} \rightarrow \frac{2}{3}$ in probability and quoting \cite[Lemma 2.4]{lawler} we have that for any $u>0$,
\begin{equation}
	\lim_{r\rightarrow \infty} \P\{V_i(u\phi_r, r) = D_r \} = \lim_{r\rightarrow \infty}\P\{U(u\phi_r,r)=D_r \}
\end{equation}
provided that one of the limits exists and equals 0 or 1.  Combining these two facts
\begin{equation}
	\P\left\{V_i\left((2/3 + \epsilon/2)\phi_r,r\right) = D_r\right\} \rightarrow 1. \nonumber
\end{equation}
In particular, (\ref{extension}) does not hold when $v = 2/3 + \epsilon$, so by the reasoning immediately preceding (\ref{extension}) there exists $c_1<\infty$ such that for all $r$:
\begin{equation}
\P\left\{V_i\left(\left(2/3 + \epsilon\right)\phi_r,r\right) \neq D_r\right\} \leq \frac{c_1}{(\log r)^2} \leq O\left(\frac{1}{\phi_r}\right).
\end{equation}
For the other bound, we follow the proof of \cite[Proposition 2.6]{lawler} to see that since  
\begin{equation}
\P\left\{V_i\left(\left(2/3 - \epsilon/2\right)\phi_r,r\right) \neq D_r\right\} \geq \frac{c_1}{(\log r)^\alpha}
\end{equation}
 for some $c_1,\alpha>1$ and all $r$ (note that the LHS converges to 1), we have:
\begin{equation}\label{est_ref}
\P\left\{V_i\left(\left(2/3 - \epsilon\right)\phi_r,r\right) = D_r\right\} \leq \E\left(1-c_2(\log r)^{-\alpha}\right)^{Y_r}
\end{equation}
where $Y_r$ is a random variable satisfying $\P(Y_r \leq \frac{1}{2}(\log r)^{\alpha + 1}) \leq O(\frac{1}{\phi_r})$ and 
$c_2>0$.  (Equation (\ref{est_ref}) is equivalent to the last equation on p. ~197 of \cite{lawler}.  The bound on the distribution of $Y_r$ follows from the last equation on p. ~198 of \cite{lawler}, once one observes that $Y_r \geq  \sum_{i=1}^J (1 - Z_i)$).  Therefore
\begin{equation}
\P\left\{ V_i\left(\left(2/3 - \epsilon\right)\phi_r,r\right)=D_r\right\} \leq e^{-c_2 \log r}[1 + O((\log r)^{1-\alpha})] + O\left(\frac{1}{\phi_r}\right) \leq O\left(\frac{1}{\phi_r}\right).
\end{equation}
\end{proof}

The next task is to relate the $\mathrm{i.i.d.}$ covering problem to the simple random walk covering problem by constructing a coupling between the two processes.  In particular, we construct the random sets $A(i,r)$ based closely on the $\mathrm{i.i.d.}$ sets so that $A(i,r)$ equals the corresponding $\mathrm{i.i.d.}$ set with high probability.  The idea behind the construction is the following.  The hitting probability on $\partial D_{2r}$ starting from any point on $\partial D_{\wp(r)}$ is nearly the same as harmonic measure.  So for $i\geq 1$ we can take $A(i,r)$ equal to a new $\mathrm{i.i.d.}$ set with probability $p$ close to 1, and equal to a random set with a biased starting position with probabilty $1-p$.  The distribution of the set $A(0,r)$ is substantially different.  However, we can couple the random walk to the $\mathrm{i.i.d.}$ set process so that $A(0,r)$ is a subset of the first $\mathrm{i.i.d.}$ set generated that contains the origin.  This allows us to bound the contribution of $A(0,r)$ from above.  A similar construction will also be used to control the second discrepancy between the coupled processes.  

We now give a precise description of the coupling.  Construct $A(0,r)$ as follows:  let $m^E = \min \{ j:0 \in E(j,r) \}$ and let $S^E$ be the random walk which was used to generate $E(m^E,r)$.  Set $t_i^E = \min \{j:S^E(j) = 0 \}$ and $t_f^E = \min \{ j: S^E(j) \in \partial D_{\wp(r)} \}$.  Define $A(0,r) = S^E[t_i^E,t_f^E] \cap D_r$.  Now, assuming that $A(0,r), \dots, A(\ell,r)$ have been constructed let $\mu_{2r}^\ell$ be the measure on $\partial D_{2r}$ which satisfies:
\begin{equation}
\mu_{2r}^\ell(x) = \P \{S(k(\ell+1,r)) = x| A(0,r), \dots, A(\ell,r) \}.
\end{equation}
That is, $\mu_{2r}^\ell$ is the hitting probability of $\partial D_{2r}$ conditioned on the subsets of $D_r$ covered by the first $\ell$ excursions. By \cite[Lemma 2.1]{lawler} we know that $|\mu_{2r}^\ell(x) - H_{2r}(x)|<\frac{c_1}{(\log r)^2}H_{2r}(x)$ where $c_1$ is independent of $\ell$ and the $A(i,r)$'s.  Thus,
\begin{equation}
\nu_{2r}^\ell = \frac{\mu_{2r}^\ell - \left(1-\frac{c_1}{(\log r)^2}\right)H_{2r}}{\frac{c_1}{(\log r)^2}}
\end{equation}
is a probability measure and we have $\mu_{2r}^\ell = (1-\frac{c_1}{(\log r)^2})H_{2r} + \frac{c_1}{(\log r)^2} \nu_{2r}^\ell$.  Define $\mathrm{i.i.d.}$ random variables $\xi_1, \xi_2, \xi_3, \dots$ so that:
\begin{eqnarray}
\P(\xi_i = 0) &=& 1-\frac{c_1}{(\log r)^2} \nonumber \\
\P(\xi_i = 1) &=& \frac{c_1}{(\log r)^2}. \nonumber
\end{eqnarray}
Then if $\xi_{\ell+1} = 0$ we set $A(\ell+1,r) = C(\ell+1,r)$.  If $\xi_{\ell+1} = 1$ and $\xi_s = 1$ for some $s \leq \ell$ we simply pick $x \in \partial D_{2r}$ according to the distribution $\nu_{2r}^\ell$, start a random walk at $x$, stop it when it hits $\partial D_{\wp(r)}$ and let $A(\ell+1,r)$ consist of the points in $D_r$ visited by the walk.  Finally, if $\xi_{\ell+1} = 1$ and $\xi_s = 0$ for all $s \leq \ell$ then we pick $x \in \partial D_{2r}$ according to $\nu_{2r}^\ell$ as before.  But now, we let $m^F = \min \{j:x \in F(j,r) \}$ and let $S^F$ be the random walk used to generate $F(m^F,r)$.  Set $t_i^F = \min \{j:S^F(j) = x \}$, $t_f^F = \min \{j:S^F(j) \in \partial D_{\wp(r)} \}$, and take $A(\ell+1,r) = S^F[t_i^F,t_f^F]$.  It is easy to see that the $A(i,r)$'s have the correct joint distributions, so this construction gives a coupling between the $\mathrm{i.i.d.}$ covering process and covering process of interest.

The motivation for this construction is as follows:  with overwhelming probability $A(i,r)$ will equal $C(i,r)$ for all $i< O(\phi_r)$ except for $i=0$ and possibly some other value, say $i_p$.  We will use the $E(i,r)$ sets to control the effect of the discrepancy at $i=0$ and the $F(i,r)$ sets to control the potential discrepancy at $i = i_p$.  Our next lemma bounds the probability that there are more than two discrepancies.

\begin{lemma}  If $u$ is a constant then $\P \left\{ \sum_{i=1}^{u\phi_r} \xi_i>1\right\} \leq O\left(\frac{1}{(\log_2 r)^{2}}\right)$.
\end{lemma}
\begin{proof}
Set $m = u\phi_r$ and $\alpha = \frac{c_1}{(\log r)^2}$.  Then
\begin{eqnarray}
	\P\left\{ \sum_{i = 1}^{u\phi_r} \xi_i>1 \right\} &=& 1 - (1-\alpha)^m - m \alpha(1-\alpha)^{m-1} \nonumber\\
	&=& 1 - [1-m \alpha+\frac{m(m-1)}{2}\alpha^2 + O(m^3 \alpha^3)] \nonumber \\
	& &-m \alpha[1-(m-1)\alpha + O(m^2 \alpha^2)] \nonumber\\
	&=& \frac{m(m-1)}{2}\alpha^2 + O(m^3 \alpha^3) \nonumber\\
	&\leq& O(m^2 \alpha^2) = O\left(\frac{1}{(\log_2 r)^2}\right).
\end{eqnarray}
\end{proof}

Next we show that if $a = \frac{64\log r \log_3 r}{\log_2 r}$ then with high probability $A(0,r) \subset V_E(a,r)$ and, if $i_p$ exists, $A(i_p,r) \subset V_F(a,r)$.  Intuitively, two discrepancies are worth no more than $2a$ excursions, which is insignificant relative to $\phi_r$.

\begin{lemma}
$\P\{m^E>a\} \leq O\left(\frac{1}{(\log_2 r)^{2}}\right)$
\end{lemma}
\begin{proof}
\begin{eqnarray}
	\P\{m^E>a\} &=& (\P \{ 0 \notin E(1,r) \})^a \nonumber \\
	&\leq& \left[1-\frac{1}{16}\left(\frac{\log(\wp(r)) - \log(2r) - c_2}{\log(\wp(r))}\right)\right]^a \nonumber 
\end{eqnarray}
by Corollary  \ref{simple}.  Thus,
\begin{eqnarray}
\P\{m^E>a\}&\leq& \left[1 - \frac{1}{32} \left(\frac{\log_2 r}{\log r}\right)\right]^a \nonumber \\
&\leq& \exp[-2 \log_3 r] \left[1 + O\left(\frac{(\log_2 r)( \log_3 r)}{\log r}\right) \right] \nonumber \\
&\leq& O\left( \frac{1}{(\log_2 r)^2}\right),
\end{eqnarray}
\end{proof}
\begin{lemma}
$\P \{m^F>a\} \leq O \left( \frac{1}{(\log_2 r)^2}\right)$
\end{lemma}
\begin{proof}
Denote by $\tau_x$ the hitting time of the point $x$ by a simple random walk.  Then
\begin{eqnarray}
	\P \{ m^F>a\} &\leq& \left[\max_{x,y \in \partial D_{2r}} \P^y \left\{ \tau_x > \zeta(\wp(r)) \right\} \right]^a \nonumber \\
&\leq& \left[ \max_{z \in \partial D_{4r}} \P^z \left\{ \tau_0>\zeta\left(\frac{1}{2}\wp(r)\right) \right\} \right]^a \nonumber \\
&\leq& \left[ 1 - \frac{1}{16}\left(\frac{\log(\frac{1}{2}\wp(r)) - \log(4n)  - c_2}{\log(\frac{1}{2}\wp(r))}\right)\right]^a \nonumber \\
&\leq& \left[1 - \frac{1}{32} \left( \frac{\log_2 r}{\log r}\right)\right]^a \nonumber \\
&\leq& O\left( \frac{1}{(\log_2 r)^2}\right)
\end{eqnarray}
by the same computation as in the previous lemma.
\end{proof}

The above lemmas provide us with sufficient control over the discrepancies between the two covering processes to establish the desired bound.  To expedite what follows, let us define 
\begin{displaymath}
Y(k,r) = \bigcup_{i=0}^k A(i,r) \;\; \mathrm{and} \;\; Z(k,r) = \bigcup_{\stackrel{1\leq i \leq k}{  \xi_i=0}} A(i,r). \nonumber
\end{displaymath}
Now compute (for sufficiently large r):
\begin{eqnarray}
	\P\left\{ Y\left(\left(2/3 + \epsilon\right)\phi_r,r\right) \neq D_r \right\} &\leq& \P\left\{Z((2/3 + \epsilon)\phi_r,r) \neq D_r \; \mathrm{and} \; \sum_{k =1}^{(\frac{2}{3} + \epsilon)\phi_r}\xi_k \leq 1 \right\} \nonumber \\
& &+ \P\left\{\sum_{k=1}^{(2/3 + \epsilon)\phi_r}\xi_k>1 \right\} \nonumber \\
&\leq& \P \{ V_C((2/3 + \epsilon/2)\phi_r,r) \neq D_r \} + O\left(\frac{1}{(\log_2 r)^2}\right) \nonumber \\
&\leq& O\left(\frac{1}{\phi_r}\right) + O\left(\frac{1}{(\log_2 r)^2}\right)  \leq O\left(\frac{1}{(\log_2 r)^2}\right).
\end{eqnarray}
The second to last inequality follows from Lemma \ref{estimate}.
For the other direction:
\begin{eqnarray} 
\P\{Y((2/3 - \epsilon)\phi_r,r) = D_r\} &\leq& \P\{V_C((2/3 - \epsilon)\phi_r,r)\cup V_E(a,r) \cup V_F(a,r) = D_r \} \nonumber \\
&& \; + \P\left\{ \sum_{k=1}^{(2/3 - \epsilon)\phi_r}\xi_k>1 \right\} + \P\{m^E>a\} + \P\{m^F>a\} \nonumber \\
&\leq& \P\{V_C((2/3 - \epsilon/2)\phi_r,r) = D_r\} + O\left(\frac{1}{(\log_2 r)^2}\right) \nonumber \\
&\leq& O\left(\frac{1}{(\log_2 r)^2}\right).
\end{eqnarray}
The second inequality holds as soon as $r$ is large enough so that 
\begin{equation}
(2/3 - \epsilon)\phi_r + 2a < (2/3 - \epsilon/2)\phi_r,
\end{equation}
 and the last inequality follows from Lemma \ref{estimate}.

\section{Asymptotic results for the Wiener sausage}
We now prove the analogous asymptotic covering result for the Wiener sausage.  Namely, if $\RR_t = \sup \{r:D(0,r) \subset \SS_t \}$ is the radius of the largest disk covered by the Wiener sausage up to time $t$ then almost surely
\begin{equation}\label{eqn34}
\limsup_{t \rightarrow \infty} \frac{(\log \RR_t)^2}{\log t \log_3 t} = \frac{1}{4}. 
\end{equation}
The style of proof is identical to the one just given.  First we give the analogues of Lemmas \ref{lem22}, \ref{lem26} and \ref{RWest}, and then show that the argument above may be duplicated with only a few minor changes to obtain (\ref{eqn34}).  Proofs of technical lemmas are collected at the end. 

Define $\zeta(r)$ to be the hitting time of $\partial D(0,r)$, $\phi_r = \frac{(\log r)^2}{\log_2 r}$ and $\wp(r) = r(\log r)^3$.  We also let $\TT_r = \inf \{t:D(0,r) \subset \SS_t \}$ be the time required for the Wiener sausage to cover $D(0,r)$, and write $\BB_t = (B^1_t,B^2_t)$ to denote the planar Brownian motion which generates $\SS_t$.  In what follows $\P^x$ denotes the probability distribution for a Brownian motion starting from $x$, if no superscript is indicated then $x=0$ is implicit.  The analogue of Lemma \ref{lem22} for Brownian motion is the following classical result:

\begin{lemma} \label{lem31}
If $r_1<r_2<r_3$ and $x \in \partial D(0,r_2)$ then we have:
\begin{equation}
\P^x\left(\zeta(r_1)<\zeta(r_3)\right) = \frac{\log r_3 - \log r_2}{\log r_3 - \log r_1}. \nonumber
\end{equation}
\end{lemma}
It is a well known fact that a Brownian motion moves \textquotedblleft twice as fast" as a random walk.  On average, a random walk requires $r^2$ steps to reach $\partial D_r$, whereas a Brownian motion reaches $\partial D(0,r)$ in time $\frac{1}{2} r^2$.  However, this factor of $\frac{1}{2}$ is insignificant for our purposes.  The following estimate follows from \cite[Theorem 6.3]{R}

\begin{lemma} \label{lem32}
	$\P \left\{ \zeta(r) \notin (r^{2-\epsilon},r^{2+\epsilon})\right\} \leq O(\exp(-r^\alpha))$ where $\alpha = \alpha(\epsilon)$.
\end{lemma}
 
To state the analogue of Lemma \ref{RWest}, we must first clarify the meaning of an excursion by the Wiener sausage. 
\begin{definition}
We shall say that an {\em excursion} from $\partial D(0,r)$ to $\partial D(0,\rho)$ by the Wiener sausage occurs between times $s<t$ if $\BB_s \in \partial D(0,r)$, $\BB_t \in \partial D(0,\rho)$ and the times $s$ and $t$ satisfy the following properties:
\begin{itemize}
\item[1.] $t = \inf \{\alpha:\alpha>s \; {\rm and} \; \BB_\alpha \in \partial D(0,\rho)\}$
\item[2.] $s = \inf \{\alpha:\alpha>\tau \; {\rm and} \; \BB_\alpha \in \partial D(0,r)\}$ where $\tau = \sup \{\alpha:\alpha<t \; {\rm and}\; \BB_\alpha \in \partial D(0,\rho)\}$ (we take the sup of the empty set to be 0).
\end{itemize}
\end{definition}
The lemma now reads:
\begin{lemma} \label{WSest}
	For any $\delta>0$ there exists $R(\delta)>0$ so that if $R\in (0,R(\delta))$ and $N_r$ represents the number of excursions between $\partial D(0,2r)$ and $\partial D(0,2R\wp(r))$ performed by the Wiener sausage before $\TT_r$, then
\begin{eqnarray}
	\P\left(N_r \leq (1-4\delta)\frac{2}{3}\phi_r\right) &\leq& O\left(\frac{(\log_2 r)^2}{\log r}\right) \nonumber \\
	\P\left(N_r \geq (1+4\delta)\frac{2}{3}\phi_r\right) &\leq& O\left(\frac{1}{r^\delta(\log r)^{3\delta}}\right).
\end{eqnarray}
\end{lemma}
The proof of this lemma constitutes the bulk of the work.  Before giving it, let us pause to verify that we may obtain (\ref{WSest}) from these three lemmas and the previous discussion.  

To prove that $\frac{1}{4}$ is an upper bound for the constant in (\ref{WSest}), we define events $A_n$ and $\tilde{A}_n$ in an analogous manner to section 2.1.  Specifically, $A_n$ is the event that $\RR_{t_{n+1}} \geq f(t_n)$ and $\tilde{A}_n$ is now the event that at least $(\frac{2}{3} - \epsilon_1)\phi_{f(t_n)}$ excursions from  $\partial D_{2f(t_n)}$ to $\partial D_{2R\wp(f(t_n))}$ occur before the Brownian path first hits $\partial D_{t_{n+1}^{1/2+\epsilon_2}}$.  Here, $R<R(\delta)$ is an arbitrary constant to be determined later.  Invoking Lemma \ref{lem31} now in place of Lemma \ref{lem22} we obtain (\ref{eqn9}), the same estimates as before yield (\ref{eqn13}) and we conclude that for any $\lambda>\frac{1}{4}$ we may choose $\epsilon_1>0$ and $\epsilon_2>0$ and $\alpha>1$ so that $\P \left\{\tilde{A}_n \; {\rm i.o.} \right\} = 0$.  To finish the argument, define $C_1$ and $D_1$ as before, and check that they have measure zero.  Lemma \ref{lem32} implies that $\P(C_1) = 0$, and Lemma \ref{WSest} gives $\P(D_1)  = 0$ once we fix $0<\delta<\frac{\epsilon_1}{2}$ and $0<R<R(\delta)$.  

For the other direction, define $A_n$ to be the event $\RR_{t_n}>f(t_n)$, and let $\tilde{A}_n$ be the event that there are at least $(2/3 + \epsilon_1)\phi_{f(t_n)}$ excursions from $\partial D_{2f(t_n)}$ to $\partial D_{2R\wp(f(t_n))}$ before $\zeta(t_n^{1/2 - \epsilon_2})$.  Again, $R<R(\delta)$ is an arbitrary constant to be determined later.  Define sigma fields $\F_n = \sigma \{\BB_t:t<\zeta(t_n^{1/2 - \epsilon_2}) \}$ and compute:
\begin{eqnarray}\label{eqn40}
	\P(\tilde{A}_n|\F_{n-1}) &\geq& \min_{\partial D(0,t_{n-1}^{1/2 - \epsilon_2})} \P^x \{\zeta(1)<\zeta(t_n^{1/2 - \epsilon_2})\}\P(\tilde{A}_n) \nonumber \\
	&\geq& \left(\frac{\log(t_n^{1/2-\epsilon_2}) - \log(t_{n-1}^{1/2-\epsilon_2})}{\log(t_n^{1/2-\epsilon_2})}\right) \P(\tilde{A}_n) \nonumber \\
	&\geq& \left(1-\frac{1}{\alpha}\right) \P(\tilde{A}_n) 
\end{eqnarray}
to obtain the analogue of (\ref{eqn16}).  Utilizing Lemma \ref{lem31} again in place of Lemma \ref{lem22} we obtain (\ref{eqn17}), and (\ref{eqn18}) follows as before.  We deduce from (\ref{eqn18}) and (\ref{eqn40}) that for any $\lambda<\frac{1}{4}$ there exist $\epsilon_1>0$ and $\epsilon_2>0$ so that $\sum_{n=1}^\infty \P(\tilde{A}_n|\F_{n-1})$ diverges a.s. and hence $\P\{ \tilde{A}_n \; {\rm i.o.} \} = 1$.  The argument is then completed by defining sets $C_2$ and $D_2$ as in section 2.2 and checking that they have measure zero.  Fixing $0<\delta<\frac{\epsilon_1}{2}$ and $0<R<R(\delta)$, this fact follows.

\subsection{Proof of Lemma \ref{WSest}}

The proof of this lemma is based largely on covering time estimates derived in \cite{DPRZ}.  Since this work deals almost exclusively with Brownian motion on a torus, we shall work primarily in that setting and then translate our results at the end.  We begin by introducing some notation, which we have chosen to be consistent with \cite{DPRZ}.  Let $\{X_t\}_{t\geq0}$ denote a Brownian motion on the torus $\T^2$ identified with the set $(-1/2,1/2]^2$.  We use the natural metric on $\T^2$ so that for any $x \in \T^2$ there is an isometry $i_x:D_{\T^2}(x,1/2) \rightarrow D(0,1/2)$.  Given a point $x \in \T^2$ define $T(x,\epsilon) = \inf \{t: d(X_t,x)\leq \epsilon \}$ to be the epsilon cover time of $x$.  For $E \subset \T^2$ define $\CC_\epsilon(E) = \sup \{ T(x,\epsilon):x \in E \}$ to be the $\epsilon$ cover time of $E$.  We shall be particularly interested in studying the number of excursions between concentric circles that occur before $\epsilon$ coverage of disks centered at the origin.  Consider $\epsilon>0$ small enough so that $2r_\epsilon<R<R(\delta)$ where $r_\epsilon = \frac{1}{|\log \epsilon|^3}$. Let $\tau^{(0)}$ be the hitting time of $\partial D_{\T^2}(0,R)$ and for $j >0$ let $\tau^{(j)}_\epsilon$ measure the sum of the durations of the $j^{th}$ excursion from $\partial D_{\T^2}(0,R)$ to $\partial D_{\T^2}(0,r_\epsilon)$, and the $(j+1)^{st}$ excursion from $\partial D_{\T^2}(0,r_\epsilon)$ to $\partial D_{\T^2}(0,R)$.  Then $N_\epsilon^\prime (R) = \max \{j:\sum_{i=0}^j \tau_\epsilon^{(i)} \leq \CC_\epsilon(D_{\T^2}(0,r_\epsilon/2)) \}$ counts the number of excursions from $\partial D_{\T^2}(0,R)$ to itself passing through $\partial D_{\T^2}(0,r_\epsilon)$ that occur before $\epsilon$ coverage of $D_{\T^2}(0,r_\epsilon/2)$.  

Our first task is to bound the probability that $N_\epsilon^\prime(R)$ differs substantially from $\frac{2}{3} \phi_\epsilon$, where $\phi_\epsilon = \frac{(\log \epsilon)^2}{\log_2 [\epsilon^{-1}]}$.  To this end, fix $1>\delta>0$ and let $N_\epsilon = \frac{2}{3}\phi_\epsilon$.  Also, define $N_\epsilon^- = \frac{2}{3}(1-4\delta)\phi_\epsilon$ and $N_\epsilon^+ = \frac{2}{3}(1+4\delta)\phi_\epsilon$.  Observe that for small enough $\epsilon$
\begin{eqnarray}
	\frac{2}{\pi}(1-\delta)(\log \epsilon)^2 &\geq& (1+\delta)\frac{N_\epsilon^-}{\pi}\log\left(\frac{R}{r_\epsilon}\right) \nonumber\\ 
	\frac{2}{\pi}(1+\delta)(\log \epsilon)^2 &\leq& (1-\delta)\frac{N_\epsilon^+}{\pi}\log\left(\frac{R}{r_\epsilon}\right) \nonumber
\end{eqnarray}
and therefore:
\begin{eqnarray} \label{eqn41}
\P(N_\epsilon^\prime(R) \leq N_\epsilon^-) &\leq& \P\left\{\CC_\epsilon(D_{\T^2}(0,r_\epsilon/2)) \leq \frac{2}{\pi}(1-\delta)(\log \epsilon)^2 \right\} \\
& & + \P\left\{\sum_{j=0}^{N_\epsilon^-} \tau^{(j)} \geq (1+\delta) \frac{N_\epsilon^-}{\pi}\log(R/r_\epsilon) \right\} \nonumber
\end{eqnarray}
\begin{eqnarray}\label{eqn42}
\P(N_\epsilon^\prime(R) \geq N_\epsilon^+) &\leq& \P\left\{\CC_\epsilon(D_{\T^2}(0,r_\epsilon/2)) \geq \frac{2}{\pi}(1+\delta)(\log \epsilon)^2 \right\} \\
& & + \P\left\{\sum_{j=0}^{N_\epsilon^+} \tau^{(j)} \leq (1-\delta) \frac{N_\epsilon^+}{\pi}\log(R/r_\epsilon) \right\}. \nonumber
\end{eqnarray}
By \cite[Lemma 2.2]{DPRZ}, we know that if $\delta>0$ is small enough and $R < R(\delta)$, then for all sufficiently small $\epsilon$ the probabilities in equations (\ref{eqn41}) and (\ref{eqn42}) involving $\tau$'s are $\leq O(\exp(-C\delta^2 N_\epsilon))$.  The following lemmas allow us to bound the other terms.
\begin{lemma} \label{lem35}
$\P\left\{\CC_\epsilon(D(0,r_\epsilon/2)) > \frac{2}{\pi}(1+\delta)(\log \epsilon)^2 \right\} \leq O(\epsilon^\delta)$
\end{lemma}
\begin{proof}
Choose a collection, $\PP$ of points in $\T^2$ so that no point in $\T^2$ is distance more than $\frac{\epsilon}{2}$ from all points in $\PP$ and $|\PP| \leq \frac{4}{\epsilon^2}$.  By Lemma 2.3 in \cite[Lemma 2.3]{DPRZ}, we know that for all $\gamma>0$ there exists $c(\gamma)<\infty$ so that for all $y>0$ and sufficiently small $\epsilon>0$ we have $\P^{x_0}(T(x,\epsilon)\geq y(\log \epsilon)^2) \leq c\epsilon^{(1-\gamma)\pi y}$ for all $x,x_0 \in \T^2$.  It follows that if $\epsilon>0$ is small enough
\begin{eqnarray}
	\P\left(\CC_\epsilon(D(0,r_\epsilon/2)) \geq \frac{2}{\pi}(1+\delta)(\log \epsilon)^2\right) &\leq& \P\left(\CC_\epsilon(\T^2) \geq \frac{2}{\pi}(1+\delta)(\log \epsilon)^2\right) \nonumber \\
&\leq& \sum_{x_0 \in \PP} \P\left\{T(x_0,\epsilon/2) \geq \frac{2}{\pi}(1+\delta)(\log \epsilon)^2 \right\} \nonumber \\
 	&\leq& \sum_{x_0 \in \PP} \P\left\{T(x_0,\epsilon/2) \geq \frac{2}{\pi}\left(1+\frac{2\delta}{3}\right)\left(\log \frac{\epsilon}{2}\right)^2 \right\} \nonumber \\
	&\leq& 4c(\gamma) \epsilon^{2(1-\gamma)(1+2\delta/3)-2} = 4c(\gamma) \epsilon^{4\delta/3 - \gamma(2+4\delta/3)}.
\end{eqnarray}
Since $\gamma>0$ is arbitrary, we obtain the desired result.
\end{proof}
\begin{lemma} \label{lem36}
$\P\left\{ \CC_\epsilon(D(0,r_\epsilon/2)) \leq (1-\delta)\frac{2}{\pi}(\log \epsilon)^2 \right\} \leq \frac{c (\log_2 1/\epsilon)^2}{\log \epsilon}$
\end{lemma}
\begin{proof}
This proof parallels the discussion in \cite[section 3]{DPRZ}.  For fixed $a<2$ and $0<\epsilon_1<R(\delta)$ define $n_k = 3ak^2\log k$, $\rho_n = n^{-25}$, $\epsilon_n = \epsilon_1 (n!)^{-3}$, and $\epsilon_{n,k} = \rho_n\epsilon_n(k!)^3$ for $1 \leq k \leq n$.  Let $R_n = \frac{1}{(\log \epsilon_n)^4}$, $S_n=[\epsilon_1 R_n,2\epsilon_1 R_n]^2$ and take $i:D_{T^2}(0,1/2) \rightarrow D(0,1/2)$ to be the natural isomorphism.  For any $x \in S_n$ denote by $\tau_n^x$ the time until $X_t$ completes $n_n$ excursions from $\partial D_{\T^2}(i^{-1}(x),\epsilon_{n,n-1})$ to $\partial D_{\T^2}(i^{-1}(x),\epsilon_{n,n})$.  Also, for $x \in S_n$ and $2 \leq k \leq n$ let $N_{n,k}^x$ denote the number of excursions of $X_t$ from $\partial D_{\T^2}(i^{-1}(x),\epsilon_{n,k-1})$ to $\partial D_{\T^2}(i^{-1}(x),\epsilon_{n,k})$ until time $\tau_n^x$.  The utility of the above construction lies in the following definition:
\begin{definition}
	A point $x \in S_n$ is said to be {\em $n$-successful} if
\begin{equation}
N_{n,2}^x = 0,  \; \;\;\; \;n_k-k \leq N_{n,k}^x \leq n_k+k \;\; \forall \;\;3,4 \dots, n-1.
\end{equation}
\end{definition}
Observe that if $x \in S_n$ is $n$-successful then $\TT(i^{-1}(x), \epsilon_{n,1}) \geq \tau_n^x$, so by studying the expected duration of excursions one can prove estimates concerning the $\epsilon_{n,1}$ cover time of $n$-successful points.  In particular, we have the following estimate from \cite[Lemma 3.2]{DPRZ}.
\begin{lemma}\label{lem38}
	There exists $C>0$ independent of $n$ so that for fixed $a<2$:
\begin{equation}
\P \left( \TT(i^{-1}(x),\epsilon_{n,1}) \leq \left(\frac{a}{\pi} - \frac{2}{\sqrt{\log n}} \right) (\log \epsilon_n)^2, x \; {\rm is}\; n{\textrm{-successful}} \right) \leq e^{-C n^2}.
\end{equation}
\end{lemma}
We shall use this lemma to control the $\epsilon_{n,1}$ covering time of $S_n$ by constructing a subset $V_n$ of $S_n$ for each $n$ and estimating the probability that $V_n$ contains an $n$-successful point.  Construct $V_n$ as follows:  partition $S_n$ into $M_n = \epsilon_1^2 R_n^2/4\epsilon_n^2$ non-overlapping squares of side length $2\epsilon_n$ and take $V_n$ to be the collection of their centers, $x_{n,j}$ for $j = 1,\dots M_n$.  Since $M_n \leq e^{o(n^2)}$ Lemma \ref{lem38} shows that for some $C_1>0$ (independent of $n$)
\begin{eqnarray}
\P\left(\CC_{\epsilon_{n,1}}(D_{\T^2}(0,R_n)) \leq (\log \epsilon_n)^2 (\frac{a}{\pi} - \frac{2}{\sqrt{\log n}} ) \right) &\leq& \P(V_n {\rm \; contains \; no \; }n{\textrm{-successful point}}) \nonumber \\
& &+ O(e^{-C_1 n^2}).\label{eqn46}
\end{eqnarray}
Now, given $\epsilon_{n+1,1}\leq \epsilon<\epsilon_{n,1}$ small enough so that $r_\epsilon/2<R_n$ we use monotonicity of cover times to write
\begin{equation}\label{eqn47}
\frac{\CC_\epsilon(D_{\T^2}(0,r_\epsilon/2))}{(\log \epsilon)^2} \geq \frac{\CC_{\epsilon_{n,1}}(D_{\T^2}(0,R_n))}{(\log \epsilon_n)^2} \frac{(\log \epsilon_n)^2}{(\log \epsilon_{n+1,1})^2}.  
\end{equation}
Observing that
\begin{equation}
\lim_{n \rightarrow \infty} \frac{(\log \epsilon_n)^2}{(\log \epsilon_{n+1,1})^2} = 1  
\end{equation}
we deduce from equations (\ref{eqn46}) and (\ref{eqn47}) that by taking $a<2$ large enough and $\epsilon>0$ small enough 
\begin{equation} \label{eqn49}
\P \left( \frac{\CC_\epsilon(D_{\T^2}(0,r_\epsilon/2))}{(\log \epsilon)^2} \leq \frac{2}{\pi}(1-\delta) \right) \leq \P(V_n {\rm \; contains \; no \;}n{\textrm{-successful point}}) + O(e^{-C_1 n^2})
\end{equation}
where $n$ is chosen to satisfy $\epsilon_{n+1,1}\leq\epsilon<\epsilon_{n,1}$.

It remains to estimate the probability that $V_n$ contains no $n$-successful point.  This computation is done in \cite{DPRZ}, but our values of $\rho_n$ and $M_n$ are different.  For $j = 1,\dots M_n$ define $Y(n,j) = 1$ if $x_{n,j}$ is $n$-successful and $Y(n,j) = 0$ otherwise.  We quote the following estimates from \cite[Lemma 3.1]{DPRZ}:
\begin{lemma} \label{lem39}
There exists $\delta_n \rightarrow 0$ such that for all $n \geq 1$,
\begin{equation}
	\bar{q}_n = \P(x \; {\rm is} \; n{\textrm{-successful}}) \geq \epsilon_n^{a + \delta_n}.
\end{equation}
For some $C_0<\infty$ and all $n$, if $|x_{n,i} - x_{n,j}| \geq 2\epsilon_{n,n}$, then
\begin{equation}
\E(Y(n,i)Y(n,j)) \leq (1 + C_0 n^{-1} \log n) \bar{q}_n^2.
\end{equation}
Further, for any $\gamma>0$ we can find $C = C(\gamma)<\infty$ so that for all $n$ and $l = l(i,j) = \max\{k\leq n: |x_{n,i} - x_{n,j}| \geq 2 \epsilon_{n,k} \} \vee 1$, 
\begin{equation}
\E(Y(n,i)Y(n,j)) \leq \bar{q}_n^2 C^{n-l}n^{39} \left(\frac{\epsilon_{n,n}}{\epsilon_{n,l+1}} \right)^{a + \gamma}.
\end{equation}
\end{lemma}
The proof of this lemma given in \cite{DPRZ} does not depend on the value of $\rho_n$, so the estimates are valid in our setting (in fact, $\bar{q}_n$ is independent of $\rho_n$).  Applying the lemma, we see that for $1\leq l \leq n-1$:
\begin{eqnarray}
	V_l &:=& (M_n\bar{q}_n)^{-2} \sum_{\stackrel{i,j = 1}{l(i,j) = l \; \mathrm{and} \; i \neq j}}^{M_n} \E\left(Y(n,i)Y(n,j)\right) \nonumber \\
	&\leq& C_1 M_n^{-1} \epsilon_{n,l+1}^2\epsilon_n^{-2}C^{n-l}n^{39} \left(\frac{\epsilon_{n,l+1}}{\epsilon_{n,n}}\right)^{-a-\gamma} \nonumber \\
	&\leq& C_1 C_2 n^{-3}(\log n)^8 C^{n-l} \left(\frac{\epsilon_{n,l+1}}{\epsilon_{n,n}}\right)^{2-a-\gamma} \nonumber \\
	&\leq& C_1 C_2 n^{-3}(\log n)^8 C^{n-l} \left(\frac{1}{(n-l-1)!^3} \right)^{2-a-\gamma}.
\end{eqnarray}
Now applying Chebyshev's inequality (see \cite[Theorem 4.3.1]{sp}) and Lemma \ref{lem39} we obtain:
\begin{eqnarray}
	\P \left( \sum_{j=1}^{M_n} Y(n,j) = 0 \right) &\leq& (M_n\bar{q}_n)^{-2} \E \left\{ \left( \sum_{i=1}^{M_n}Y(n,i)\right)^2\right\} -1 \nonumber \\
	&\leq& (M_n \bar{q}_n)^{-1} + C_0 n^{-1}\log n + \sum_{l = 1}^{n-1} V_l \nonumber \\
	&\leq&  C_1 n^8(\log n)^8 \epsilon_n^{2-a-\delta_n} + C_0 n^{-1}\log n \nonumber \\
& & + C_3 n^{-3}(\log n)^8 \sum_{j=1}^\infty C^j \left(\frac{1}{j!}\right)^{6-3a-3\gamma} 
\end{eqnarray}
Fixing $\gamma>0$ so that $2-a-\gamma>0$ it follows that the probability that no $x \in V_n$ is $n$-successful is $O(n^{-1}\log n)$.  If $\epsilon_{n+1,1}<\epsilon<\epsilon_{n,1}$ then $\log (1/\epsilon) \sim 3n \log n$ and $\log_2 (1/\epsilon) \sim \log n$ so by (\ref{eqn49}) we obtain
\begin{equation}
	\P\left(\CC_\epsilon(D_{\T^2}(0,r_\epsilon/2)) \leq (1-\delta)\frac{2}{\pi}(\log \epsilon)^2 \right) \leq O\left(\frac{(\log_2 \epsilon)^2}{\log \epsilon}\right)
\end{equation}
as desired.
\end{proof}

Equations (\ref{eqn41}) and (\ref{eqn42}) together with Lemmas \ref{lem35} and \ref{lem36} now imply that for $0<R < R(\delta)$:
\begin{eqnarray}
	\P\left\{N_\epsilon^\prime(R) \leq N_\epsilon^- \right\} &\leq& O\left(\frac{(\log_2 \epsilon)^2}{\log \epsilon} \right)\label{eqn56} \\
	\P\left\{N_\epsilon^\prime(R) \geq N_\epsilon^+ \right\} &\leq& O\left( \epsilon^\delta \right)\label{eqn57}
\end{eqnarray}
(recall that the probabilities in (\ref{eqn41}) and (\ref{eqn42}) involving $\tau$'s are $O(\exp(-C\delta^2N_\epsilon))$ so they are insignificant).  

The next task is to extend these results to planar Brownian motion.  Let $N_\epsilon(R)$ be the analogue of $N_\epsilon^\prime(R)$ for planar Brownian motion, that is $N_\epsilon(R)$ is the number of $\BB_t$ excursions between $\partial D(0,R)$ and itself passing through $\partial D(0,r_\epsilon)$ that occur prior to $\epsilon$ coverage of $D(0,r_\epsilon/2)$.  It follows from the discussion in \cite[p. 18]{DPRZ} that if $R>0$ is taken small enough there is a function $h_{N_\epsilon, r_\epsilon}$ for all $\epsilon<\epsilon(R)$ such that
\begin{eqnarray}
	\P(N_\epsilon(R) \leq N_\epsilon^-) &\leq& ||h_{N_\epsilon^-,r_\epsilon}||_\infty \P(N_\epsilon^\prime(R)\leq N_\epsilon^-) \label{eqn58} \\
	\P(N_\epsilon(R) \geq N_\epsilon^+) &\leq& ||h_{N_\epsilon^+,r_\epsilon}||_\infty \P(N_\epsilon^\prime(R)\geq N_\epsilon^+)  \label{eqn59}
\end{eqnarray}
 with $||h_{N_\epsilon^+,r_\epsilon}||_\infty, ||h_{N_\epsilon^-,r_\epsilon}||_\infty \rightarrow 1$ as $\epsilon \rightarrow 0$.  Thus, we deduce the analogues of (\ref{eqn56}) and (\ref{eqn57}) for planar Brownian motion:
\begin{eqnarray}
	\P\left\{N_\epsilon(R) \leq N_\epsilon^- \right\} &\leq& O\left(\frac{(\log_2 \epsilon)^2}{\log \epsilon} \right) \label{eqn60} \\
	\P\left\{N_\epsilon(R) \geq N_\epsilon^+ \right\} &\leq& O\left( \epsilon^\delta \right)\label{eqn61}
\end{eqnarray}
provided that $R < R(\delta)$ and $R(\delta)$ is taken small enough for (\ref{eqn58}) and (\ref{eqn59}) to apply.  Brownian scaling implies that (\ref{eqn60}) and (\ref{eqn61}) still hold if $N_\epsilon(R)$ is redefined to count the number of excursions between $\partial D(0,\frac{1}{\wp(\epsilon)})$ and $\partial D(0,R/\epsilon)$ necessary to cover $D(0,\frac{1}{2\wp(\epsilon)})$ by the Wiener sausage.  Now fix $\epsilon_r$ so that $r = \frac{1}{2\wp(\epsilon_r)}$.  Fix $R<R(\delta)$ and $\gamma>0$ so that $(1+\gamma)R<R(\delta)$.  Define $N_r^{\pm}= \frac{2}{3}(1 \pm 4\delta)\phi_r$ and let $N_r(R)$ be the random variable counting the number of excursions between $\partial D(0,2r)$ and $\partial D(0,2R \wp(r))$ that occur prior to the $\TT_r$.  Observing that $\log \epsilon_r \sim \log r$ and $\epsilon_r \sim \frac{1}{2r (\log r)^3}$ we see that for large $r$:
\begin{eqnarray}
\P(N_r(R) \leq N_r^-) &\leq& \P(N_{\epsilon_r}((1+\gamma)R) \leq N_{\epsilon_r}^-) \leq O\left(\frac{(\log_2 \epsilon_r)^2}{\log \epsilon_r}\right) = O\left(\frac{(\log_2 r)^2}{\log r}\right)  \\
\P(N_r(R) \geq N_r^+) &\leq& \P(N_{\epsilon_r}((1-\gamma)R) \geq N_{\epsilon_r}^+) \leq O\left(\epsilon_r^\delta\right) = O\left(\frac{1}{r^\delta(\log r)^{3\delta}}\right)  
\end{eqnarray}
and we obtain Lemma \ref{WSest}.

\noindent{\bf Acknowledgements.} 
We are grateful to Olivier Daviaud, Alan Hammond and Manjunath Krishnapur for useful comments.  Some of this research was conducted while the first author was visiting the University of Washington and the second author was visiting Microsoft Research.  We thank them for their hospitality.

\bigskip
\noindent
\sc \bigskip \noindent J. Ben Hough, Department
of Mathematics, U.C.\ Berkeley, CA 94720,
USA.\\{\tt jbhough@math.berkeley.edu},
{\tt www.math.berkeley.edu/\~{}jbhough

\sc \bigskip \noindent Yuval Peres, Departments of
Statistics and Mathematics, U.C.\ Berkeley, CA 94720,
USA. \\{\tt peres@stat.berkeley.edu}, {\tt
stat-www.berkeley.edu/\~{}peres}

\end{document}

%% file: macros.tex
\def\Z{{\mathbb Z}}        
\def\T{{\mathbb T}}		

\def\F{{\mathcal F}}
\def\E{{\mathbb E}}
\def\P{{\mathbb P}}

\def\1{{\mathbf 1}}

\def \BB{{\mathcal B}}
\def \SS{{\mathcal S}}
\def \RR{{\mathcal R}}
\def \TT{{\mathcal T}}

\def \CC{{\mathcal C}}

\def \PP{{\mathcal P}}

